\documentclass[letterpaper, 10 pt, conference]{ieeeconf}
\IEEEoverridecommandlockouts
\usepackage{amssymb,amsmath,color}
\usepackage{graphicx}
\usepackage{epsfig}
\usepackage{algorithm,algorithmic,xspace}

\renewcommand{\natural}{{\mathbb{N}}}
\newcommand{\naturalzero}{\natural_0}

\newcommand{\real}{{\mathbb{R}}}
\newcommand{\subscr}[2]{{#1}_{\textup{#2}}}
\newcommand{\union}{\cup}
\newcommand{\intersection}{\ensuremath{\operatorname{\cap}}}

\newcommand{\map}[3]{#1: #2 \rightarrow #3}

\newcommand{\setdef}[2]{\{#1 \; | \; #2\}}

\newcommand{\vers}{\operatorname{vers}}
\newcommand{\sign}{\operatorname{sign}}
\newcommand{\fl}[1]{\left\lfloor #1\right\rfloor}
\newcommand{\ceil}[1]{\left\lceil #1\right\rceil}
\newcommand{\diam}{\operatorname{diam}}

\newcommand{\Enorm}[1]{\|#1\|_{2}}

\newcommand{\umax}{\subscr{r}{ctr}}
\newcommand{\rcmm}{\subscr{r}{cmm}}

\newcommand{\GG}{\mathcal{G}}
\newcommand{\cball}[2]{B(#2,#1)}
\newcommand{\ccube}[2]{C(#2,#1)}

\newcommand{\until}[1]{\{1,\dots,#1\}}

\newcommand{\NN}{\mathcal{N}}

\newcommand{\RR}{\mathcal{R}}

\newcommand{\finite}{\mathbb{F}}

\newcommand{\zeroRd}{{0}_{\real^d}}

\renewcommand{\zeroRd}{{0}_d}

\newcommand{\subj}{\text{subj. to}}

\newcommand{\minimize}{\text{minimize}}
\newcommand{\half}{\frac{1}{2}}

\newcommand{\normI}[1]{\|#1\|_{\infty}}

\newcommand{\MEB}{\operatorname{MEB}}
\newcommand{\MBC}{\operatorname{MBC}}
\newcommand{\MEO}{\operatorname{MEO}}
\newcommand{\MOC}{\operatorname{MOC}}

\newcommand{\agent}{A}
\newcommand{\setofagents}{\mathcal{A}}
\newcommand{\true}{\textup{\texttt{true}}\xspace}
\newcommand{\false}{\textup{\texttt{false}}\xspace}

\newcommand{\network}[1][]{\Sigma_{\textup{#1}}}
\renewcommand{\network}[1][]{\mathcal{S}_{\textup{#1}}}
\newcommand{\supind}[2]{{#1}^{[#2]}}
\newcommand{\FCC}[1][]{\mathcal{CC}_{\textup{#1}}}
\newcommand{\task}[1][]{\mathcal{T}_{\textup{#1}}}

\newcommand{\nll}{\textup{\texttt{null}}\xspace}
\newcommand{\ctrl}{\textup{ctl}}
\newcommand{\msg}{\textup{msg}}

\newcommand{\stf}{\textup{stf}}

\newcommand{\myboolean}{\textup{\{true, false\}}}

\newcommand{\rendezvous}{rndzvs}

\newcommand{\TC}{\operatorname{TC}}

\newcommand{\MTR}{\mathcal{MTR}}
\newcommand{\MTrendset}{\subscr{\RR}{MT}}

\newcommand{\MEBradius}{\operatorname{MBR}}

\newtheorem{theorem}{Theorem}[section]

\newtheorem{definition}[theorem]{Definition}
\newtheorem{lemma}[theorem]{Lemma}
\newtheorem{remark}[theorem]{Remark}

\def\QEDclosed{\mbox{\rule[0pt]{1.3ex}{1.3ex}}} 

\def\QED{\QEDclosed} 

\newcommand\oprocendsymbol{\hbox{$\square$}}
\newcommand\oprocend{\relax\ifmmode\else\unskip\hfill\fi\oprocendsymbol}
\def\eqoprocend{\tag*{$\square$}}

\renewcommand{\theenumi}{(\roman{enumi})}

 \renewcommand{\baselinestretch}{0.975}

\begin{document}

\title{Distributed consensus on enclosing shapes\\
  and minimum time rendezvous\thanks{This material is based upon work
    supported in part by ARO MURI Award W911NF-05-1-0219. The authors would
    like to thank Sonia Mart{\'\i}nez, Jorge Cort\'es, and Emilio Frazzoli
    for numerous discussions on robotic networks.
}}

\author{Giuseppe Notarstefano\thanks{Giuseppe Notarstefano is with the the
    Department of Information Engineering, Universit\`a di Padova, Via
    Gradenigo 6/b, 35131 Padova, Italy, \texttt{notarste@dei.unipd.it}}
  \and
  Francesco Bullo
\thanks{Francesco Bullo is with the
    Center for Control, Dynamical Systems and Computation, University of
    California at Santa Barbara, Santa Barbara, CA 93106,
    \texttt{bullo@engineering.ucsb.edu}}
}

\maketitle

\begin{abstract}
  In this paper we introduce the notion of optimization
  under control and communication constraint in a robotic
  network. Starting from a general setup, we focus our
  attention on the problem of achieving rendezvous in
  minimum time for a network of first order agents with
  bounded inputs and limited range communication.  We
  propose two dynamic control and communication laws.  These
  laws are based on consensus algorithms for distributed
  computation of the minimal enclosing ball and orthotope of
  a set of points. We prove that these control laws converge
  to the optimal solution of the centralized problem (i.e.,
  when no communication constrains are enforced) as the
  bound on the control input goes to zero.  Moreover, we
  give a bound for the time complexity of one of the two
  laws.
\end{abstract}

\section{Introduction}
The interesting aspect of motion coordination consists in combining
together problems from control and communication theory. The main
difficulty deals with integrating the sensing, computing, communication and
control aspects of problems involving groups of mobile agents. A well known
problem in control theory is optimal control. Roughly speaking, it consists
in finding a feedback law that minimizes some cost functional under some
inputs and dynamics constraint. In this paper we introduce the notion of
optimal control and communication for a network of robotic agents. We want
to study how to solve an optimization problem, in presence of both the
usual motion constraints and the communication ones. In particular this
paper is a preliminary contribution towards what might be loosely referred
to as ``distributed geometric optimization.'' In fact many optimization
problems for robotic networks can be shown to be equivalent to the
computation of geometric shapes. While in a centralized setting the
solution is usually simple, the problem becomes very complicated when it
must be solved in a distributed way. Distributed computation over network
has been largely studied for fixed topologies; e.g., see~\cite{NAL:97}.

In this paper we point our attention on the well known rendezvous
coordination task and look for solutions that solve such task in
minimum time. We look for distributed solutions in networks of
mobile agents with first order dynamics, bounded inputs and
limited-range communication.

The ``multi-agent rendezvous'' problem and a ``circumcenter algorithm''
have been introduced by Ando and coworkers in~\cite{HA-YO-IS-MY:99}.  The
algorithm proposed in~\cite{HA-YO-IS-MY:99} has been extended to various
synchronous and asynchronous stop-and-go strategies
in~\cite{JL-ASM-BDOA:03}. A related algorithm, in which connectivity
constraints are not imposed, is proposed in \cite{ZL-MB-BF:04a}. In
\cite{JC-SM-FB:04h} the class of ``circumcenter algorithms'' has been
studied in networks of agents whose state space is $\real^d$, for arbitrary
$d$, and with communication topology characterized by proximity graphs
spatially distributed over the disk graph. In \cite{SM-FB-JC-EF:05mn-tmp}
the (time and communication) complexity of this and other algorithms has
been studied. All these coordination schemes are memoryless (static
feedback). In this paper we want to explore dynamic control and
communication laws in order to approximate the optimal solution of the
minimum time rendezvous. In particular the control and communication laws
is based on reaching consensus on some logic variables and at the same time
moving toward the current estimation. A similar approach was used in
\cite{WR-RWB-TWM:04} where the agents try to reach a consensus on a set of
variables called coordination variables.

Studying the minimum time rendezvous problem in the centralized
setting we show that, depending on the norm used to bound the
control input, the optimal solution consists of moving toward the
center of the minimal enclosing ball (bound on $L_2$ norm) and
toward the center of the minimal enclosing orthotope (bound on the
infinity norm) of the points located at the initial position of the
agents.

Our main result is the design of a control and communication law
based on a consensus algorithm for the distributed computation of
the minimal enclosing ball and the minimal enclosing orthotope of a
set of points. We prove the correctness of the two consensus
algorithms and provide a bound on the time of convergence for the
orthotope case. Then we prove that the control and communication law
that combines the consensus with the motion law converges to the
optimal solution as the control bound goes to zero. Moreover, for
the problem with input bounded by the infinity norm (corresponding
to the computation of the minimal enclosing orthotope), we prove
that the control and communication law is a constant factor
approximation of the centralized optimal solution.

In Section~\ref{sec:prelim} we introduce a formal model of robotic
network inspired by the one introduced in
\cite{SM-FB-JC-EF:05mn-tmp}. Moreover, we define the optimal control
and communication problem. In Section~\ref{sec:mintime_rend} we
characterize the solution of minimum time rendezvous in a
centralized setting. In Section~\ref{sec:consensus-MEB-MEO} we
define the \emph{FloodMEB} and \emph{FloodMEO} algorithms for the
distributed computation of minimal enclosing ball and orthotope,
prove their correctness and give bounds on time complexity.
Section~\ref{sec:contr-comm-laws} contains the control and
communication laws based on the consensus algorithms described in
Section~\ref{sec:consensus-MEB-MEO}. Finally in
Section~\ref{sec:simulations} and Section~\ref{sec:conclusions} we
show simulations and draw the conclusions with future perspectives.


\section{Preliminary developments}
\label{sec:prelim} In this section we recall the concepts of network
of robotic agents, coordination tasks and complexity measures, and
introduce the notion of optimization under motion and communication
constraints.

\subsection{Notation}
We let $\natural$, $\naturalzero$, and $\real_+$ denote the natural
numbers, the non-negative integer numbers, and the positive real numbers,
respectively. We let $ \prod_{i\in\until{n}}S_i$ denote the Cartesian
product of sets $S_1,\dots, S_n$.  For $p\in\real$, we let $\fl{p}$ and
$\ceil{p}$ denote the floor and ceil of $p$.
For $r\in\real_+$ and $p\in\real^d$, we let $\cball{r}{p}$ denote the
closed ball centered at $p$ with radius $r$, i.e.,
$\cball{r}{p}=\setdef{q\in\real^d}{\Enorm{p-q}\leq{r}}$ and $\ccube{r}{p}$
denote the closed hypercube centered at $p$ with sides of length $r$ and
parallel to the coordinate axes, i.e.,
$\ccube{r}{p}=\setdef{q\in\real^d}{\normI{p-q}\leq{r}}$.

For $\map{f,g}{\natural}{\real}$, we say that $f \in O(g)$
(respectively, $f \in \Omega(g)$) if there exist $n_0 \in
\natural$ and $k \in \real_+$ such that $|f(n)| \le k |g(n)|$ for
all $n \ge n_0$ (respectively, $|f(n)| \ge k |g(n)|$ for all $n
\ge n_0$).  If $f \in O(g)$ and $f \in \Omega(g)$, then we use the
notation $f \in \Theta (g)$.

Next, we briefly review some useful proximity graphs.
%
%
%
%
Given $\rcmm\in\real_+$, the \emph{disk graph} $\subscr{\GG}{disk}(\rcmm)$,
respectively \emph{cube graph} $\subscr{\GG}{cube}(\rcmm)$, is the state
dependent graph on $\real^d$ defined by the following statement: for any
pointset $\{\supind{p}{1},\dots,\supind{p}{n}\}\subset\real^d $, the pair
$(i,j)$ is an edge in
$\subscr{\GG}{disk}(\rcmm)\cdot(\{\supind{p}{1},\dots,\supind{p}{n}\})$,
respectively
$\subscr{\GG}{cube}(\rcmm)\cdot(\{\supind{p}{1},\dots,\supind{p}{n}\})$, if
and only if $i\neq j$ and
\begin{equation*}
  \Enorm{\supind{p}{i} - \supind{p}{j} } \le \rcmm
  \enspace \iff \enspace  \supind{p}{i} - \supind{p}{j}\in
  \cball{\rcmm}{\zeroRd},
\end{equation*}
respectively
\begin{equation*}
  \|  \supind{p}{i} - \supind{p}{j} \|_{\infty} \le \rcmm
  \enspace \iff \enspace  \supind{p}{i} - \supind{p}{j}\in \ccube{\rcmm}{\zeroRd}.
\end{equation*}
Another useful graph is the complete graph $\subscr{\GG}{cmpl}$, i.e., the
graph with edges between any pair of nodes.

Finally, given a graph $\GG$ (even not state dependent), we denote
with $\text{dist}_\GG(i, j)$ the topological distance between $i$
and $j$, i.e., the minimum number of agents to go from $i$ to $j$ in
the graph $\GG$. We define $\diam_{\GG}$, the \emph{diameter} of
$\GG$, to be the maximum topological distance, $\text{dist}_\GG(i,
j)$, for all $(i, j)$.

\subsection{Modeling a network of robotic agents}
We describe a (uniform) network of robotic agents using the formal model
introduced in \cite{SM-FB-JC-EF:05mn-tmp} modified for the discrete time
case. The network is modeled as a tuple $(I,\setofagents,\subscr{E}{cmm})$.
$I=\until{n}$ is the \emph{set of unique identifiers (UIDs)};
$\setofagents=\{\supind{\agent}{i}\}_{i\in{I}}= \{(X,U,X_0,f)\}_{i\in{I}}$
is called the \emph{set of physical agents} and is a set of control systems
consisting of a differentiable manifold $X$ (state space), a compact subset
$U$ of $\real^m$ (input space), a subset $X_0$ of $X$ (set of allowable
initial states) and a (sufficiently smooth) map $\map{f}{X \times U}{X}$
describing the dynamics of $i$th agent;
$\map{\subscr{E}{cmm}}{X^n}{I\times{I}}$ is called the \emph{communication
  edge map}.

The robotic network evolves according to a discrete-time
communication and motion model.
\begin{definition}[Control and communication law]
  \label{dfn:dynamic-feedback-control-communication}%
  Let $\network$ be a robotic network.  A \emph{(uniform, synchronous, dynamic)
    control and communication law} $\FCC$ for $\network$ consists of the
  sets:
  \begin{enumerate}

  \item $L$, a set containing the $\nll$ element, called the
    \emph{communication language}; elements of $L$ are called
    \emph{messages};

  \item $W$, set of values of some \emph{logic
      variables} $\supind{w}{i}$, $i\in{I}$;

  \item $W_{0} \subseteq W$,
    subsets of \emph{allowable initial values};
  \end{enumerate}
  and of the maps:
  \begin{enumerate}
  \item $\map{ \msg }%
    {X \times W \times I}
    {L}$, \emph{message-generation function};

  \item $\map{\stf}%
    {W \times L^n } {W}$,
    called \emph{state-transition function};

  \item $\map{\ctrl}%
    {X \times
      W \times L^n } {U}$, called
    \emph{control function}. \oprocend
\end{enumerate}
\end{definition}

Roughly speaking this definition has the following meaning: for
all $i\in{I}$, to the $i$th physical agent corresponds a logic
process, labeled $i$, that performs the following actions. First,
at each communication round the $i$th logic process sends to each
of its neighbors in the communication graph a message (possibly
the $\nll$ message) computed by applying the message-generation
function to the current values of $\supind{x}{i}$ and
$\supind{w}{i}$.  After a negligible period of time, the $i$th
logic process resets the value of its logic variables
$\supind{w}{i}$ by applying the state-transition function to the
current value of $\supind{w}{i}$, and to the messages received at
time $t$. Between communication instants, the motion of the $i$th
agent is determined by applying the control function to the
current value of $\supind{x}{i}$, and the current value of
$\supind{w}{i}$. This idea is formalized as follows.

\begin{definition}[Evolution of a robotic network%
] \label{dfn:evolution}%
  Let $\network$ be a robotic network and $\FCC$ be a control and
  communication law for $\network$.  The \emph{evolution} of
  $(\network,\FCC)$ from initial conditions $\supind{x}{i}_{0}\in
  X_0$ and $\supind{w}{i}_{0}\in W_0$, $i\in{I}$,
  is the set of curves
  $\map{x^{[i]}}{\natural}{X}$ and
  $\map{\supind{w}{i}}{\natural}{W}$, $i\in{I}$,
  satisfying
  \begin{align*}
    x^{[i]}(t+1) &= f\big( x^{[i]}(t), \,
    \ctrl(x^{[i]}(t), \supind{w}{i}(t+1), y^{[i]}(t) ) \big) ,
  \end{align*}
  where, for $i\in{I}$,
  \begin{align*}
   \supind{w}{i}(t+1) =
    \stf(\supind{w}{i}(t), y^{[i]}(t)) \,
    ,
  \end{align*}
  with the conventions that $x^{[i]}(t_0)=\supind{x}{i}_{0}$ and
  $\supind{w}{i}(t_{0})=\supind{w}{i}_{0}$, $i\in{I}$. Here, the
  function $\map{\supind{y}{i}}{\natural}{L^n}$ (describing the
  messages received by agent $i$) has components
  \begin{equation*}
      \supind{y}{i}_j(t) =
    \begin{cases}
      \msg(x^{[j]}(t),w^{[j]}(t),i),
      & \text{if} \enspace (i,j)\in
      \subscr{E}{cmm},
      \\
      \nll, & \text{otherwise}. \eqoprocend
  \end{cases}
\end{equation*}
\end{definition}

In the paper we consider the following network. Each agent $i$
occupies a location $\supind{p}{i} \in \real^d$, $d \in \natural$,
and moves according to the first order discrete-time integrator
\begin{equation}\label{eq:agent-dynamics}
  \supind{p}{i}(t+1) = \supind{p}{i}(t) + \supind{u}{i}(t).
\end{equation}
The communication edge map can be either the one arising according
to the \emph{disk graph}, $\subscr{E}{disk}$, or the one according
to the \emph{cube graph}, $\subscr{E}{cube}$. Each control
$\supind{u}{i}$ takes values in a bounded subset of $\real^d$,
that can be either $\cball{\umax}{0}$ or $\ccube{\umax}{0}$, i.e.,
$\Enorm{\supind{u}{i}} \le \umax$ or $\normI{\supind{u}{i}} \le
\umax$. Notice that, in general, the type of communication edge
map and the type of control bound are not related. Finally the
control and communication law will be defined depending on the
coordination task.

\subsection{Coordination tasks and time complexity}
We are ready to define the notion of task and of task achievement
by a robotic network.

\begin{definition}[Coordination task] \label{dfn:task}
  Let $\network$ be a robotic network.
 A \emph{(static) coordination task} for $\network$ is a map $\map{\task}{
      X^n}{\myboolean}$.
  Additionally, let $\FCC$ a control and communication law for $\network$.
  The law $\FCC$ \emph{achieves} the task $\task$ if, for all initial conditions
    $\supind{x}{i}_{0}\in X_0$ and $\supind{w}{i}_0\in
    W_{0}$, $i\in I$, the corresponding network evolution
    $t\mapsto(x(t), w(t))$ has the property that there exists $T\in\natural$
    such that $\task(x(t)) = \true$ for all $t\geq T$.  \oprocend
\end{definition}

We are finally ready to define the notion of time complexity as
the minimum number of communication rounds needed by the agents to
achieve the task $\task$ with $\FCC$.

\begin{definition}[Time complexity] \label{dfn:time-complex}
  Let $\network$ be a robotic network and let $\task$ be a coordination
  task for $\network$.  Let $\FCC$ be a control and communication law for
  $\network$ compatible with $\task$.  The \emph{time complexity to achieve
    $\task$ with $\FCC$ from $x_0 \in
    X^n_0$}
  is
$$    \TC(\task,\FCC,x_0) = \inf \,  \{T\in\natural\;|\;
    \task(x(t))=\true\,,\; \forall \; t \ge T\}$$
  where $t\mapsto(x(t),w(t))$ is the evolution of $(\network,\FCC)$ from
  the initial condition $(x_0,w_0)$.

  The \emph{time complexity to achieve $\task$ with $\FCC$}, $\TC (\task,\FCC)$, is
  the maximum $\TC(\task,\FCC,x_0)$ over all initial conditions $x_0$.
\end{definition}

\subsection{Optimal control and communication in robotic networks}
Having defined a coordination task for a robotic network, we can
ask whether such task can be accomplished minimizing some cost
functional. In what follows we will introduce the notion of
\emph{optimal control and communication problem} and of
\emph{optimal control and communication law} as solution of the
problem.
\begin{definition}[Optimal control and communication]
\label{dfn:optimal-control-communication-problem} Given a task
$\task$ and a cost functional $J(u(\cdot), x(T), T)$, an optimal
control and communication problem is the
following:\\

\parbox{150mm}{
  $\minimize_{u(\cdot), x(0), x(T), T} ~~ J(u(\cdot), x(T), T)\\[1.5ex]
  J(u(\cdot), x(T), T) = \sum_{\tau = 0}^{T}( l(x(\tau), u(\tau)) +
  g(x(T)),$ \\[1.5ex]
  $\subj$
  \begin{enumerate}
  \item $(x(\cdot),u(\cdot))$ is an input-state trajectory of $\setofagents$,\\
    $\setofagents=\{\supind{\agent}{i}\}_{i\in{I}}
    $;
  \item $i$ and $j$ can communicate if and only if\\
    $(i,j)\in \subscr{E}{cmm}(\supind{x}{1}(t),\dots,\supind{x}{n}(t))$;
  \item $\task(x(t)) = \true$ for all $t\ge T$, $T \in \natural$.
  \end{enumerate}}\\

\noindent where $\map{l}{X^n  \times
  U^n }{\real}$ is a sufficiently smooth and
nonnegative-valued function, called \emph{stage cost}, and
$\map{g}{X^n }{\real}$ has the same properties plus $g(x)=0$ for all
$x\in X^n $ such that $\task(x) = \true$ (for an admissible $\FCC$).
\oprocend
\end{definition}

We say that a control and communication law $\FCC$ is optimal with
respect to the coordination task $\task$ and the cost functional
$J$, if it solves the above optimal control and communication
problem.


We call $\FCC$ a \emph{centralized} optimal control and
communication law if it solves the optimization problem for a
network of robotic agents that communicate according to the
complete graph, i.e., the communication edge map is
$\subscr{E}{cmpl}$.

\begin{remark}
The centralized solution of an optimal control and communication
problem is the classical solution of the optimal control problem
for the whole network system without communication constraints.
\oprocend
\end{remark}

\section{Centralized minimum time rendezvous}
\label{sec:mintime_rend} In this section we study the rendezvous
problem for a robotic network of first order agents with
communication edge map $\subscr{E}{disk}$ or $\subscr{E}{cube}$
and look for a control and communication law that solves the
problem in minimum time.

More formally, let $\network=(I,\setofagents,\subscr{E}{cmm})$ be
a uniform robotic network. The \emph{(exact) rendezvous task}
$\map{\task[\rendezvous]}{X^n}{\myboolean}$ for $\network$ is the
static task defined by
\begin{equation*}
    \task[\rendezvous] (x) =
    \begin{cases}
      \true, & \text{if} \enspace
      \supind{x}{i}=\supind{x}{j}, \;\\
      & \forall (i,j)\in \subscr{E}{cmm}(x),\\
      \false, & \text{otherwise} .
    \end{cases}
  \end{equation*}
  for $x=(\supind{x}{1},\dots,\supind{x}{n})$.

  Thus, given the uniform network
  $\network=(I,\setofagents,\subscr{E}{cmm})$, the \emph{minimum time
    rendezvous} problem for first order
  agents with limited-range communication and bounded control input is the following:\\

\parbox{150mm}{
$\minimize \sum_{\tau = 0}^{T} 1,$\\[1.5ex]
$u(\cdot), p(T)$\\[1.5ex]
$\subj$
\begin{enumerate}
\item $(p(\cdot),u(\cdot))$ is an input-state trajectory of
$\setofagents$,\\
$\setofagents=\{\supind{\agent}{i}\}_{i\in{I}}=
\{(\real^d, U,\real^d,f)\}_{i\in{I}}, ~p(0) = p_0$;
\item $i$ and $j$ can communicate if and only if\\
        $(i,j)\in \subscr{E}{cmm}(\supind{p}{1}(t),\dots,\supind{p}{n}(t))$;
\item $\task[\rendezvous] (\supind{p}{1},\dots,\supind{p}{n}) = \true$ for all $t\ge T$,
$T \in \natural$.
\end{enumerate}
}

\noindent Here $U$ is either $\cball{\umax}{0}$ or
$\ccube{\umax}{0}$, $f(\supind{p}{i}(t), \supind{u}{i}(t)) =
\supind{p}{i}(t) + \supind{u}{i}(t)$ and the communication edge
map $\subscr{E}{cmm}$ is either $\subscr{E}{disk}$ or
$\subscr{E}{cube}$.

We refer to the minimum time rendezvous problem with communication
edge map $\subscr{E}{cmm}$ and input set $U$ as
$\MTR(\subscr{E}{cmm}, U)$.

Next, we provide some preliminary results for the centralized setting of
the above problem, that is, for $\MTR(\subscr{E}{cmpl}, U)$. Let
$\MEB(\supind{p}{1} \cdots \supind{p}{n})$ and $\MEO(\supind{p}{1} \cdots
\supind{p}{n})$ the minimal enclosing ball and orthotope of points
$(\supind{p}{1} \cdots \supind{p}{n})$, and let $\MBC(\supind{p}{1} \cdots
\supind{p}{n})$ and $\MOC(\supind{p}{1} \cdots \supind{p}{n})$ the centers
of $\MEB(\supind{p}{1} \cdots \supind{p}{n})$ and $\MEO(\supind{p}{1}
\cdots \supind{p}{n})$ respectively.  We present the following theorem
omitting the proof based on geometric arguments because of space
constraints.

\begin{theorem}
For all $\umax \in \real_+$, $\supind{p}{i}_0 \in \real^d$, $i \in
\until{n}$ the solution of $\MTR(\subscr{E}{cmpl},U)$,
$U=\cball{\umax}{0}$ or $U=\ccube{\umax}{0}$, is not unique (the
problem is not normal). If

\noindent $\supind{u}{i} \in \cball{\umax}{0}$, $i\in\until{n}$,
then
\newcounter{saveenum}

\begin{enumerate}

\item $p(T) =\subscr{p}{\rendezvous, disk} = \MBC(\supind{p}{1}(0), \dots, \supind{p}{n}(0))$,
\begin{equation*}
\begin{split}
\supind{u}{i}(t) =& \min\{\umax,
\Enorm{\subscr{p}{\rendezvous}-\supind{p}{i}(t)}\} ~~~~~~~~~~~~~~~~~~~~~~\\
& \cdot \vers{(\subscr{p}{\rendezvous}-\supind{p}{i}(t))}, \quad i\in\until{n},
\end{split}
\end{equation*}
is a solution of $\MTR(\subscr{E}{cmpl},\cball{\umax}{0})$;
\item if $\forall i \in \until{n}, ~ \Enorm{\supind{p}{i} - \MBC(\supind{p}{1} \cdots
\supind{p}{n})} \le \umax$, then the solution of
$\MTR(\subscr{E}{cmpl},\cball{\umax}{0})$ is given by $p(T)=
\subscr{p}{\rendezvous, disk}$, $\subscr{p}{\rendezvous, disk}\in
\intersection_{i\in\until{n}} \cball{\umax}{\supind{p}{i}}$, and
$\supind{u}{i}(t) =\subscr{p}{\rendezvous, disk} -
\supind{p}{i}(t)$.
\setcounter{saveenum}{\value{enumi}}
\end{enumerate}

\noindent Alternatively, if $\supind{u}{i} \in \ccube{\umax}{0}$,
$i\in\until{n}$, then
\begin{enumerate}
\setcounter{enumi}{\value{saveenum}}
\item $p(T) =\subscr{p}{\rendezvous, cube}$, $\subscr{p}{\rendezvous,
cube}\in \MTrendset$,
 $\supind{u_{a}}{i}(t) = \min\{\umax, |\subscr{p}{\rendezvous,a}-\supind{p_a}{i}(t)|\}
\sign(\subscr{p}{\rendezvous,a}-\supind{p_a}{i}(t))$,
$i\in\until{n}$, $a\in\until{d}$ is a solution of
$\MTR(\subscr{E}{cmpl},\ccube{\umax}{0})$, where
\begin{equation*}
\begin{split}
\MTrendset =& \prod_{a}\left[\MOC(\supind{p}{1}(0), \dots,
\supind{p}{n}(0)) - \half(\subscr{l}{max} -
\subscr{l}{a}),\right.\\
& \left. \MOC(\supind{p}{1}(0), \dots, \supind{p}{n}(0)) + \half
(\subscr{l}{max} - \subscr{l}{a})\right],
\end{split}
\end{equation*}

$\subscr{l}{max}$ is the largest side of $\MEO(\supind{p}{1}(0),
\dots, \supind{p}{n}(0))$ and $\subscr{l}{a}$ is the side in direction
$a$;

\item if $\forall i \in \until{n} ~ \normI{\supind{p}{i} - \MOC(\supind{p}{1} \cdots
\supind{p}{n})} \le \umax$ the solution of
$\MTR(\subscr{E}{cmpl},\ccube{\umax}{0})$ is given by $p(T)=
\subscr{p}{\rendezvous}$, $\subscr{p}{\rendezvous}\in
\intersection_{i\in\until{n}} \ccube{\umax}{\supind{p}{i}}$, and
$\supind{u}{i}(t) =\subscr{p}{\rendezvous} - \supind{p}{i}(t)$. \oprocend
\end{enumerate}
\end{theorem}

\section{Distributed consensus on minimal enclosing ball and orthotope}
\label{sec:consensus-MEB-MEO} In the previous section we have shown
that minimal enclosing shapes play a key role in the solution of
minimum time rendezvous. In fact if the agents could know the center
of such shapes (ball or orthotope) the solution of minimum time
rendezvous would be just a control law that drives each agent to
this point. Therefore, in this section, we want to explore two
consensus algorithms to compute the minimal enclosing ball and the
minimal enclosing orthotope of a set of given points in $\real^d$ in
a distributed way.


Here is an informal description of what we shall refer to as the
\emph{FloodMEB} algorithm:
\begin{quote}
  \emph{[Informal description]} Each agent initializes the minimal
  enclosing ball to its initial position, then, at each
  communication round, performs the following tasks: (i) it acquires
  from its neighbors (a message represented by) the coordinates of
  the minimum set of points describing the boundary of their minimal
  enclosing ball and the coordinates of their initial position;
  (ii) it computes the minimal enclosing ball of the
  point set comprised of its and its neighbors' set of points and
  its initial position (that it maintains in memory);
  (iii) it updates its logic variables and message as in (i).
\end{quote}


Before describing the algorithm more formally, we need to introduce
some notation and state some properties of the minimal enclosing
ball. Given a set of $m$ points $\{q_1,\dots,q_m\}\subset\real^d$ in
generic positions, we denote with $\subscr{\MEB}{bndry}(\{q_1,
\dots, q_m\})$ the minimum set of points on the boundary of
$\MEB(\{q_1, \dots, q_m\})$ that uniquely identify such boundary.
When the points are in generic position, we let
$\subscr{\MEB}{bndry}(\{q_1, \dots, q_m\})$ denote \emph{a} minimum
set of points on the boundary of $\MEB(\{q_1, \dots, q_m\})$ that
identify such boundary. 
Moreover, let $\map{\MEBradius}{\finite{(\real^d)}}{\real}$ the
function that associates to a set of points the radius of the
minimal enclosing ball of such points.

\begin{lemma}[MEB properties]
\label{lmm:MEB-monotonicity} Let $Q_n$ 
a set of $n$ points. The following statements hold.
\begin{enumerate}
\item there exists a subset $Q_d \subset Q_n$ of $d+1$
elements such that $\MEB(Q_d) = \MEB(Q_n)$;
\item for all $Q_{n_1}, Q_{n_2} \subset Q_n$ with $Q_{n_1} \subset Q_{n_2}$, then $\MEBradius(Q_{n_1}) \le
\MEBradius(Q_{n_2})$;
\item if $\MEBradius(Q_{n_1}) = \MEBradius(Q_{n_2})$, then $\MBC(Q_{n_1}) =
\MBC(Q_{n_2})$;
\item the number of possible values of $\MEBradius(Q_{n_1})$, for all
  $Q_{n_1} \subset Q_n$, is finite. \oprocend
\end{enumerate}
\end{lemma}

\begin{remark}
An important implication of Lemma~\ref{lmm:MEB-monotonicity}(i) is
that $\subscr{\MEB}{bndry}(\{q_1, \dots, q_n\})$ has at most $d+1$
points, then the number of packets in the message sent and stored
by each agent is at most $d+1$ and does not depend on $n$.
\oprocend
\end{remark}

The algorithm is described formally in the following table.

\begin{table}[htbp]
  \centering
  \noindent{
    \normalsize
    \framebox[.9999\linewidth]{\noindent\parbox{.9999\linewidth-2\fboxsep}{%
        \noindent\begin{tabular}{ll}
          \textbf{Name:}      & \parbox[t]{.65\linewidth}{\emph{FloodMEB} algorithm.}\\[2.5ex]
          \textbf{Goal:}      &
          \parbox[t]{.65\linewidth}{Solve the problem of computing minimal enclosing ball of a set of points.}
          \\[1.5ex]
          \textbf{Logic state:} & $\supind{w}{i} =
          (\subscr{\supind{P}{i}}{bndry}, \supind{p}{i}_0);
          $\\[1.5ex]
          \textbf{Msg function:} & $\msg(\supind{x}{i}, \supind{w}{i}, i) =
          \supind{w}{i}$;
          \\[1.5ex]
          \textbf{Initialization:} & $\subscr{\supind{P}{i}}{bndry}(0) = \{\supind{p}{i}(0)\}$,\\
          & $\supind{p}{i}_0(0) = \supind{p}{i}(0)$.
          \\[2ex]
          \hline
      \end{tabular}\\[.5ex]

      For $i\in\until{n}$, agent $i$ executes at each time $t\in\natural$:\\[-2ex]
      \begin{algorithmic}[1]
        \STATE acquire $\supind{w}{j}(t)$, $j \in \NN(i)$
        \\[.5ex]
        \STATE compute\\
        $\subscr{\supind{P}{i}}{bndry}(t+1)=
        \subscr{\MEB}{bndry}(W_{\NN(i)}(t))$,\\
        $W_{\NN(i)}(t)=\setdef{\supind{w}{j}(t)}{j \in \NN(i) \union i}$
        \\[.5ex]
        \STATE update
        $\supind{w}{i}(t+1) =
        (\subscr{\supind{P}{i}}{bndry}(t+1), \supind{p}{i}_0(t))$
        \\[.5ex]
      \end{algorithmic} }}}
\end{table}

\begin{remark}
For the algorithm to converge it is important that each agent keeps
in memory the coordinates of its initial position and thus computes
the minimal enclosing ball on the points received from its neighbors
together with the point located in its initial position. In fact a
point on the boundary on the minimal enclosing ball of $n_{1}$
points is not ensured to be on the boundary of the ball of $n_{2}
\le n_{1}$ points. This means that the coordinates of the agents on
the boundary could be taken out from the logic variables during the
first iterations. This does not happen, for example, for the minimal
enclosing orthotope. The result is a simplified consensus algorithm.
\oprocend
\end{remark}

We are now ready to prove the algorithm's correctness.



\begin{theorem}[Correct $\MEB$ computation]
  \label{thm:Correct-MEB-computation}
  Let $\network$ be a robotic network such that the agents can communicate
  according to some communication edge map $\subscr{E}{cmm}$. For any
  $\FCC$ such that the graph remains connected along the evolution, the
  \emph{FloodMEB} algorithm achieves consensus on minimal enclosing ball.
  \oprocend
\end{theorem}

 \begin{proof}
   In order to prove correctness of the algorithm, observe, first of all,
   that each law at every node converges in a finite number of steps. In
   fact, using Lemma~\ref{lmm:MEB-monotonicity}, each sequence corresponds
   to a ball whose radius is monotone nondecreasing, upper bounded and can
   assume a finite number of values. Then we proceed by contradiction to
   prove that all the laws converge to the same ball (same radius and
   center) and that it is exactly the minimal enclosing ball of the $n$
   points.  Suppose that the algorithm converges on different balls
   (different radius or different center) for different agents. Then there
   must exist two agents that are neighbors and have different logic
   variables (corresponding to different balls). But this means that, at the
   following time instant, they have to compute the minimal enclosing ball
   of a larger set of points, then either one of them will take the value of
   the other or both of them will change their value and take a common one.
   Iterating this argument we obtain that all the agents must converge to a
   common value. Now, the ball at each node contains, by construction, the
   initial position of that node. Since the ball is the same for each node,
   it contains all the initial positions, then it is the minimal enclosing
   ball of the initial positions.
 \end{proof}

\begin{remark}
If we admit that the agents have different priority, the initial
positions of the agents can be shared by all the agents in time of
order $\Theta(n^2)$. The algorithm is the following. Each agent
sends the position of the agent with higher (or equivalently
lower) priority that he has in memory. Each position takes time
$\Theta(n)$ to spread in the network, therefore the total time
complexity is $\Theta(n^2)$.
Even if we did not provide any bound for the time complexity of
\emph{FloodMEB} algorithm, however simulations suggest that it
should be of order $\Theta(n)$. Moreover, while the algorithm for
sharing the initial position needs to store a number of packets of
order $\Theta(n)$, the \emph{FloodMEB} algorithm needs to store
only $d+2$ packets.\oprocend
\end{remark}

Here is an informal description of what we shall refer to as the
\emph{FloodMEO} algorithm:
\begin{quote}
  \emph{[Informal description]} Each agent initializes the minimal
  enclosing orthotope to its initial position, then, at each communication
  round, performs the following tasks: (i) it acquires from its neighbors a
  message represented by the coordinates of their current minimal enclosing
  orthotope; (ii) it computes the minimal of its and its neighbors'
  enclosing orthotopes; (iii) it stores as new message the coordinates of
  the minimal enclosing orthotope computed at the previous step.
\end{quote}

\noindent A more formal description of the algorithm is given in the
  following table.

\begin{table}[htbp]
  \centering
  \noindent{
    \normalsize
    \framebox[.9999\linewidth]{\noindent\parbox{.9999\linewidth-2\fboxsep}{%
        \noindent\begin{tabular}{ll}
          \textbf{Name:}      & \parbox[t]{.65\linewidth}{\emph{FloodMEO} algorithm.}\\[2.5ex]
          \textbf{Goal:}      &
          \parbox[t]{.65\linewidth}{Solve the problem of computing minimal enclosing orthotope of a set of points.}
          \\[1.5ex]
          \textbf{Logic state:} & $\supind{w}{i} = \{ \supind{w}{i}_{a}
          \}_{a \in \until{d}}$\\
          & \phantom{$\supind{w}{i}$} $= \{ (\supind{p}{i}_{min,a},
          \supind{p}{i}_{max,a})
          \}_{a \in \until{d}}$\\[1.5ex]
          \textbf{Msg function:} & $\msg(\supind{x}{i}, \supind{w}{i}, i) =
          \supind{w}{i}$
          \\[1.5ex]
          \textbf{Initialization:} & $\supind{p}{i}_{min,a}(0) =
          \supind{p}{i}_{a}(0)$,\\
          & $\supind{p}{i}_{max,a}(0) = \supind{p}{i}_{a}(0)$
          \\[2ex]
          \hline
      \end{tabular}\\[.5ex]

      For $i\in\until{n}$, agent $i$ executes at each time $t\in\natural$:\\[-2ex]
      \begin{algorithmic}[1]
        \STATE acquire $\supind{w}{j}$, $j \in \NN(i)$
        \\[.5ex]
        \STATE compute $\forall a \in\until{d}$\\
        $\supind{p}{i}_{min,a}(t+1) = \min_{j\in\NN(i)\union\{i\}}
        \{\supind{p}{j}_{min,a}(t)\}$\\
        $\supind{p}{i}_{max,a}(t+1) = \max_{j\in\NN(i)\union\{i\}}
        \{\supind{p}{j}_{max,a}(t)\}$
        \\[.5ex]
        \STATE update $\forall a \in\until{d}$\\ $ \supind{w}{i}_{a}(t+1)
           =  \left(\supind{p}{i}_{min,a}(t+1),
          \supind{p}{i}_{max,a}(t+1) \right)$
        \\[.5ex]
      \end{algorithmic} }}}
\end{table}
In the following theorem we prove the correctness of this
algorithm, together with the fact that it reaches consensus in
minimum time.

\begin{theorem}[Correct $\MEO$ computation]
\label{thm:Correct-MEO-computation} Let $\network$ be a robotic
network such that the agents can communicate according to some
communication edge map $\subscr{E}{cmm}$. For any $\FCC$ such that
the graph remains connected along the evolution, then the
\emph{FloodMEO} algorithm achieves consensus on minimal enclosing
orthotope. Moreover, it achieves consensus in minimum number of
communication rounds given by
\begin{equation*}
\label{eq:min-rounds-orthotope} \subscr{T}{\emph{FloodMEO}}
=\max_{a\in\until{d}} \max_{i\in\until{n}} \{ \text{dist}_\GG(i,
\subscr{i}{min,a}), \text{dist}_\GG(i, \subscr{i}{max,a}) \},
\end{equation*}
where $\subscr{i}{min,a}$ and $\subscr{i}{max,a}$ are the agents
that characterize the boundary of the orthotope in direction $a$
and minimize the topological distance from $i$.

The time complexity of the algorithm is of order $\Theta(n)$.  \oprocend
\end{theorem}

 \begin{proof}
   In order to prove the correctness and the time complexity of the
   algorithm described in Table~\ref{tab:MEO-algorithm}, we need to prove
   that it is equivalent to $2d$ \emph{FloodMax} algorithms for leader
   election (two for each direction) running simultaneously.  Once we have
   proven that, the results on correctness and time complexity follow from
   Chapter~4 in \cite{NAL:97}.

   The algorithm is clearly a set of \emph{FloodMax} algorithms for leader
   election. In fact the boundary of the orthotope in each direction $a$ is
   given by the coordinates of the points on such boundary which are
   characterized by the property of having the maximum and minimum value of
   the $a$th coordinate respectively.

   In order to prove that the exact number of communication rounds needed is
   $\subscr{T}{\emph{FloodMEO}}$, simply observe that it is exactly the
   minimum time for all the leaders to propagate their information through
   all the network. Hence this is the minimum time for every possible
   consensus algorithm to converge. But this is exactly the time taken by
   $2d$ \emph{FloodMax} algorithms running simultaneously and therefore the
   time taken by \emph{FloodMEO}.
 \end{proof}

In the following lemma we give, for both \emph{FloodMEB} and
\emph{FloodMEO} algorithms, a bound on the time needed by each
agent to decide that the algorithm has reached consensus.

\begin{lemma}[Termination condition]
Consider a network $\network$, where the \emph{FloodMEB}
(\emph{FloodMEO}) algorithm is running. Each agent can decide that
the algorithm has reached consensus if the value of its $\MEB$
($\MEO$) has not changed after $\diam_{\GG}$ communication rounds.
\oprocend
\end{lemma}

 \begin{proof}
   In order to prove the claim we proceed by contradiction.  Suppose that
   after $\diam_{\GG}$ communication rounds the $\MEB$ ($\MEO$) of agent $i$
   (for some $i\in \until{n}$) has not changed and the algorithm has not
   converged yet. Then there will exist a $T > \diam_{\GG}$ such that the
   $\MEB$ ($\MEO$) of agent $i$ will change to a new value. But this means
   that the new value, stored $T$ rounds before by some other agent $j$,
   took a number of communication rounds greater than $\diam_{\GG}$ to
   arrive from $j$ to $i$ and this contradicts the definition of diameter of
   $\GG$.
 \end{proof}

\section{Constant factor approximation of minimum time rendezvous control and communication law}
\label{sec:contr-comm-laws} The centralize solution for minimum time
rendezvous and the consensus algorithms studied in the previous
section suggest a dynamic control and communication law that plays a
key role in the minimum time rendezvous problem.

Here is an informal description of what we shall refer to as the
\emph{move-toward-$\MBC$ ($\MOC$) control and communication law},
$\FCC[MEB]$ ($\FCC[MEO]$):
\begin{quote}
  Each agent initializes its logic variables to its initial position, then,
  at each communication round, performs the following tasks: (i) it
  acquires from its neighbors a message given by their logic variables and
  positions; (ii) it runs, as state transition function, the
  \emph{FloodMEB(MEO)} algorithm;  (iii) it moves toward the center of the
  current ball (orthotope) while maintaining connectivity.
\end{quote}

Next, we formally define the law as follows. First we assume that
each agent operates according to the standard message-generation
function, that is
$\msg(\supind{x}{i},\supind{w}{i},i)=(\supind{x}{i},\supind{w}{i})$.
Second, before the \emph{FloodMEB} (or \emph{FloodMEO}) algorithm
reach consensus, connectivity is maintained by restricting the
allowable motion of each agent in some appropriate manner. The exact
algorithm can be found for example in \cite{SM-FB-JC-EF:05mn-tmp}.

%
The state transition function implements the \emph{FloodMEB} and
\emph{FloodMEO} algorithms respectively, with logic variables as defined in
the two tables above.

Define the control function
$\map{\ctrl}{\real^d\times\real^d\times{L^n}} {\real^d}$ for each
agent $i\in\until{n}$ by:
\begin{equation*}
\begin{split}
  \label{eq:MEB-control-law}
  \ctrl(\supind{p}{i}, \supind{w}{i}, \supind{y}{i}) =
  \max\{\lambda_* \cdot  (p_{\rendezvous}(\supind{w}{i},\supind{y}{i}) - \supind{p}{i}), \umax \}\,\\
  \cdot  \vers(p_{\rendezvous}(\supind{w}{i},\supind{y}{i})- \supind{p}{i}) ,
\end{split}
\end{equation*}
with
\begin{equation}
\label{eq:p-rendezvous}
p_{\rendezvous}(\supind{w}{i},\supind{y}{i}) =
\MBC(\supind{w}{i},\supind{y}{i})
\end{equation}
and $\lambda_*$ is chosen in order to maintain connectivity until
consensus is reached.
%

In a network with communication edge map
$\subscr{E}{cmm}=\subscr{E}{cube}$ the procedure described above
is applied separately in every direction $a\in\until{d}$.


The correctness of the two control and communication laws is
proven in the following lemma.

\begin{lemma}[Correctness of ${\FCC[MEB]}$ and ${\FCC[MEO]}$]
On the network $\network$ with communication edge map
$\subscr{E}{ball}$ or $\subscr{E}{cube}$ and bound on the $i$th
control input $\supind{u}{i} \in \cball{\umax}{0}$ or
$\supind{u}{i} \in \ccube{\umax}{0}$, the control and
communication laws $\FCC[MEB]$ and $\FCC[MEO]$ achieve rendezvous
at $\MBC(\supind{p}{1}(0),\dots,\supind{p}{n}(0))$ and
$\MOC(\supind{p}{1}(0),\dots,\supind{p}{n}(0))$ respectively.
\oprocend
\end{lemma}

 \begin{proof}
   By the connectivity arguments done before and by
   Theorem~\ref{thm:Correct-MEB-computation} and
   Theorem~\ref{thm:Correct-MEO-computation} we know that there exists
   $\overline{T}\in \natural$ such that for $t=\overline{T}$ the network is
   connected and all the agents have reached consensus on
   $\MBC(\supind{p}{1}(0),\dots,\supind{p}{n}(0))$ (or
   $\MOC(\supind{p}{1}(0),\dots,\supind{p}{n}(0))$). Since this instant all
   the agents can move toward the same point (at maximum speed) without
   enforcing connectivity constraint anymore. Thus, they can converge to the
   rendezvous point which is exactly
   $\MBC(\supind{p}{1}(0),\dots,\supind{p}{n}(0))$ (or
   $\MOC(\supind{p}{1}(0),\dots,\supind{p}{n}(0))$).
 \end{proof}


\subsection{Time complexity of $\FCC[MEB]$ and $\FCC[MEO]$}
In the previous lemma we have proven that the control and
communication laws $\FCC[MEB]$ and $\FCC[MEO]$ achieve consensus.
Now we ask how fast these laws are depending on the control bound
$\umax$ and the number of agents.

\begin{theorem}
\label{thm:time-complexity} For $\rcmm\in\real_+$, $d\in
\natural$, consider the network $S$ with communication edge map
either $\subscr{E}{disk}$ or $\subscr{E}{cube}$. The following
statements hold:
\begin{enumerate}
\item for $\supind{u}{i}\in\cball{\umax}{0}$, $i\in\until{n}$,
the control and communication law $\FCC[MEB]$ asymptotically
converges to the minimum time rendezvous centralized solution
$\MTR(\subscr{E}{cmpl},\cball{\umax}{0})$ as $\umax \rightarrow
0^+$ (for all fixed $n$).

\item for $\supind{u}{i}\in\ccube{\umax}{0}$, $i\in\until{n}$,
the control and communication law $\FCC[MEO]$ converges to the
minimum time rendezvous centralized solution
$\MTR(\subscr{E}{cmpl},\ccube{\umax}{0})$ for $\umax \rightarrow
0^+$ (for all fixed $n$). Moreover, it is a constant factor
approximation of $\MTR(\subscr{E}{cmpl},\ccube{\umax}{0})$, i.e.,
$\TC (\task[\rendezvous],\FCC[MEO]) \in \Theta(\frac{n}{\umax})$
for $\umax \rightarrow 0^+$ and $n \rightarrow +\infty$. \oprocend
\end{enumerate}
\end{theorem}


Before proving the theorem let us state a useful lemma.

\begin{lemma}[\cite{KF-BG:04}]
  For all pointsets $P_1 \subset P$, we have $\MBC(P_1) \in \MEB(P)$ and
  $\MOC(P_1) \in \MEO(P)$. \oprocend
\end{lemma}
%
%


 \begin{proof}[Theorem~\ref{thm:time-complexity}]
   The line of proof of the two statements is the same, hence we prove only
   the second one which is stronger due to the stronger result on the time
   complexity of the \emph{FloodMEO} algorithm.  Using the previous lemma we
   know that for all $t \in \naturalzero$, then
   $\supind{p}{i}_{\rendezvous}(t)\in\MEO(\supind{p}{1}(0),\dots,\supind{p}{n}(0))$,
   where $\supind{p}{i}_{\rendezvous}$ is defined as in
   \eqref{eq:p-rendezvous}. This implies that, once the consensus is
   reached, the time to rendezvous is upper bounded by the time of the
   centralized solution. Hence the following bound on the time of
   convergence of $\FCC[MEO]$, $T_{\MEO}$, holds:
   \begin{equation}
     \label{eq:bound-on-T_MEO}
     T_{\MEO} \leq
     \ceil{\diam(\supind{p}{1}(0),\dots,\supind{p}{n}(0)) \cdot
       \frac{1}{\umax}} + \subscr{T}{\emph{FloodMEO}}.
   \end{equation}
   The first statement is proven by observing that
   $\subscr{T}{\emph{FloodMEO}}$ does not depend on $\umax$, therefore, as
   $\umax \rightarrow 0^+$, $T_{\MEO}$ converges to the optimal value of the
   centralized case.

   In order to prove the second statement, observe that
   $\diam(\supind{p}{1}(0),\dots,\supind{p}{n}(0))\le(n-1)\rcmm$ and
   $\subscr{T}{\emph{FloodMEO}} \in \Theta(n)$. The result follows by
   substituting these bounds in \eqref{eq:bound-on-T_MEO}.
 \end{proof}

\begin{remark}
The previous theorem confirms the intuitive idea that, if the
communication is much faster than the motion ($\umax$ small), then
the optimal solution in the distributed case converges to the one
of the centralized case.\oprocend
\end{remark}

\subsection{Distributed minimum time rendezvous in one dimension}
In one dimension (all the agents spread on a line), we can find a
condition on $\umax$ ensuring that the move-toward-$\MBC$
algorithm is the solution of
$\MTR(\subscr{E}{disk},\cball{\umax}{0})$.

\begin{theorem}
  For $d=1$, let $\subscr{i}{max}$ and $\subscr{i}{min}$ the agents in the
  network $\network$ with the maximum and minimum positions.  If $\umax <
  \frac{1}{4}\rcmm$ and both $\max_{j\in
    \NN(\subscr{i}{max})}\{\supind{p}{\subscr{i}{max}}(0) -
  \supind{p}{j}(0) \} \ge 2 \umax$ and $\max_{j\in
    \NN(\subscr{i}{min})}\{\supind{p}{j}(0) -
  \supind{p}{\subscr{i}{min}}(0)\} \ge 2 \umax$, then the control and
  communication law $\FCC[MEB]$ solves the task $\task[\rendezvous]$ in
  minimum time which is exactly the time of the centralized solution, that
  is, $T^* = \ceil{\frac{\diam\left(\supind{p}{1}, \dots,
        \supind{p}{n}\right)}{ \umax}}$.\oprocend
\end{theorem}

 \begin{proof}
   Consider the input sequence of the centralized solution for the agent
   $\subscr{i}{min}$ (and equivalently for $\subscr{i}{max}$). It is
   $\supind{u}{\subscr{i}{min}}(t) = \umax$ for all $t < T^*-1$ and
   $\supind{u}{\subscr{i}{min}}(T^*-1) = \MBC(\supind{p}{1}(0), \dots,
   \supind{p}{n}(0)) - \supind{p}{\subscr{i}{min}}(T^*-1)$. Since the
   rendezvous time is bounded by the time that $\subscr{i}{max}$ and
   $\subscr{i}{min}$ take to reach $\MBC(\supind{p}{1}(0), \dots,
   \supind{p}{n}(0))$, we need to prove that, as long as the consensus on
   the minimal enclosing ball is not reached, then
   $\supind{u}{\subscr{i}{min}}(t) = -\supind{u}{\subscr{i}{max}}(t) =
   \umax$. Due to the symmetry of the problem we will give the proof only
   for $\subscr{i}{min}$. It can be easily shown that for all $t \ge 1$ such
   that $\supind{p}{\subscr{i}{min}}_{max}(t) \not=
   \supind{p}{\subscr{i}{max}}(0)$ (consensus is not reached), the following
   holds:
   \begin{equation*}
     \supind{p}{\subscr{i}{min}}_{max}(t+1) >
     \supind{p}{\subscr{i}{min}}_{max}(t-1) + \rcmm.
   \end{equation*}
   It follows:
   \begin{equation*}
     \supind{\MBC}{\subscr{i}{min}}(t+1) = \half(
     \supind{p}{\subscr{i}{min}}_{max}(t+1) +
     \supind{p}{\subscr{i}{min}}(0)),
   \end{equation*}
   then
   \begin{equation*}
     \begin{split}
       \supind{\MBC}{\subscr{i}{min}}(t+1) & > \half(
       \supind{p}{\subscr{i}{min}}_{max}(t-1) + \rcmm +
       \supind{p}{\subscr{i}{min}}(0))\\
       & = \supind{\MBC}{\subscr{i}{min}}(t-1) + \half \rcmm\\
       & \ge \supind{\MBC}{\subscr{i}{min}}(t) - \umax + \half \rcmm.
     \end{split}
   \end{equation*}
   This leads to
   \begin{equation*}
     -\umax + \half \rcmm > \umax,
     \enspace\text{and}\enspace
     \umax  < \frac{1}{4} \rcmm.
   \end{equation*}
   The other two assumptions ensure the condition for $t=0$.
 \end{proof}

\subsection{Simulations} \label{sec:simulations} 
\begin{figure}[htbp]
  \centering
  \includegraphics[width=.8\linewidth]{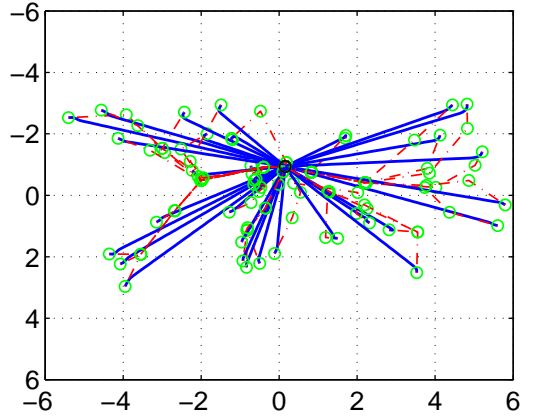}%
  \caption{Evolution of the network (in filled blue) according to
    $\FCC[MEB]$ with evolution of \emph{FloodMEB} (green circles connected
    by dashed red line)}\label{fig:run-MEB-rend}
\end{figure}

In order to
illustrate the performance of our rendezvous algorithms, we
implemented the move-toward-$\MBC$ algorithm, based on the
\emph{FloodMEB} consensus algorithm. We implemented it in the plane,
$d=2$, over the disk graph. The simulation run is illustrated in
Figure~\ref{fig:run-MEB-rend}.  The $32$ agents have a bound on the
control inputs $\umax = 0.1$, and a communication radius $\rcmm=3$.
The initial positions of the agents were randomly generated over the
rectangle
$[-6,6]\times[-3,3]$.%

The \emph{FloodMEB} law converges in five steps, while the
rendezvous is achieved at $T = 58$. As it clearly appears in the
figure, once the consensus on the minimal enclosing ball is
reached, all the agents move toward the center.

\section{Conclusions}
\label{sec:conclusions} We have presented some simple algorithms on
how to compute optimal enclosing shapes for pointsets via
distributed computation. These algorithms are then used to provide
efficient solutions to distributed rendezvous problems for
synchronous robotic networks. For future work we envision
characterizing the time complexity of the \emph{FloodMEB} algorithm
and, in turn, of the move-toward-$\MBC$ control and communication
law.


{\small
\bibliographystyle{IEEEtran}
\bibliography{brevalias,FB,New,Main}
}

\end{document}